# Cyclic behavior of masonry barrel vaults strengthened through Composite Reinforced Mortar, considering the role of the connection with the abutments

Ingrid Boem[a]* and Natalino Gattesco[a]

[a]Department of Engineering and Architecture, University of Trieste, p.le Europa 1, 34127, Trieste, Italy
*corresponding author – boem@dicar.units.it

## Abstract

The original results of experimental investigations concerning the reduction of the seismic vulnerability of masonry vaults through Composite Reinforced Mortar (CRM) are presented in the paper. The transversal performances of full scale, masonry barrel vault samples (i.e. running bond pattern vaults and brick *in folio* vaults) carrying their own weight and reinforced at the extrados or at the intrados are investigated through cyclic tests. The connection of the reinforced vaults with the masonry abutment is particularly considered and discussed, due to the importance of this detail to avoid uplift and slip at the skewback sections and better exploit the CRM benefits. The experimental behavior of the vault samples is described and analyzed in terms of global behavior, crack pattern and load-displacement curves, also in comparison with the unreinforced configurations. The CRM allowed significant improvements in both vaults strength and displacement capacities, as contrasted the cracks opening and foster the spreading of cracks all along the vault. Considerable performances emerged also in terms of dissipative capacities, with a mean damping value of about 13% in the post-cracking configuration. The vault-to-wall connection, based on embedded steel and composite bars, resulted essential for avoiding dangerous vault sliding at skewback sections, that may reduce the vault performances and also cause the fall from the support.

# 1    Introduction

Solid brick vaults were frequently employed in historic masonry buildings, with multiple variations depending on the historical period and the geographical location [1]. Many factors affect their structural behavior, such as geometry, construction technology, mass distribution, boundary conditions, over time modifications, materials deterioration. They can be either load-bearing structures (with the vault supporting a backfill or parallel spandrels on which the upper floor insists), or secondary elements carrying their own weight only (constituting the ceiling of a room).

The paper focuses on these latter structures: they are, typically, single layer, slender vaults with bricks arranged in running bond pattern or laid flatly (brick "*in folio*" vaults [2]). Their behavior during a seismic event is typically governed by bending [3] and it is strongly influenced by the interaction with the entire structure, through the abutments [4]. In particular, in slender vaults supporting its own weight and laying on bulky, fixed abutments, the collapse mechanism affects the vault only (Figure 1a); differently, when relative displacements between the abutments occurred (e.g. presence of a relatively slender wall, with no steel ties), also the vertical walls are involved (Figure 1b). In the case of fixed abutments, the activation of the collapse mechanism, concerning the formation of four structural hinges, in alternate position intrados-extrados [5]-[7], is abrupt (fragile failure), since the masonry tensile resistance is negligible and the stabilizing vertical load in the vault is very low (ceiling vaults). Thus, the lateral performances of thin vaults carrying their own weight are very weak during seismic events [8]-[9] - Figure 2.  In case of moving supports, it is firstly necessary to use ties or other effective strategies to contrast the relative displacements between the abutments.

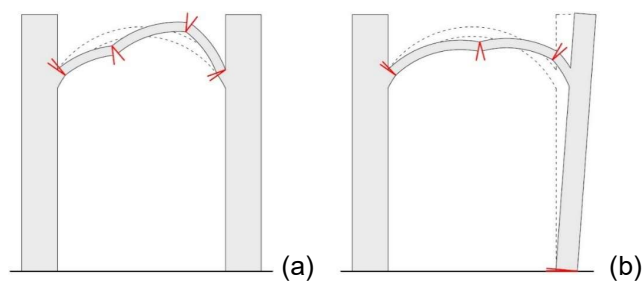

Figure 1. Typical four-hinges collapse mechanism of the masonry arch system (a) on fixed abutments or (b) with moving support.

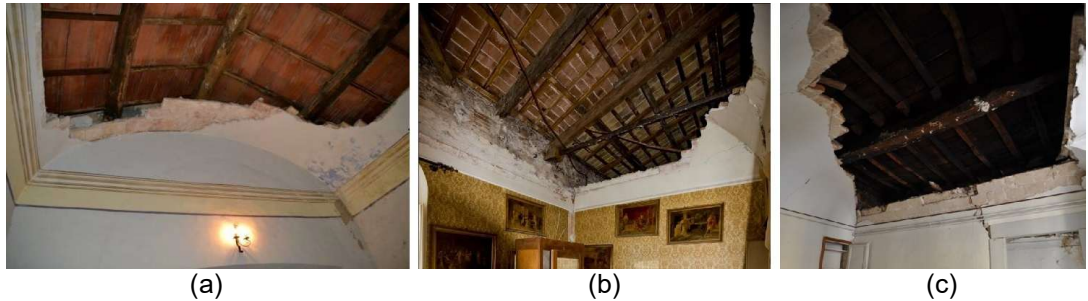

|     (a)     |     (b)     |     (c)     |

Figure 2. The failure of masonry vaults occurred during the 2016 earthquake in Central Italy, (a) in the St. Egidio's Church in Guardea [2] and (b)-(c) in some residential buildings (ph. courtesy of prof. A. Borri, University of Perugia).

Several strengthening solutions were contrived and applied for the reinforcement of arch and vaulted elements; strength and weaknesses were already broadly analyzed and discussed by many authors [10]-[15]. Typically, the concept of these strategies is to introduce, at the extrados or/and at the intrados, tensile resistant elements, in the form of externally bonded layers or structural repointing. The aim is to improve the masonry bending performances, avoiding the opening of the localized structural hinges in favor of a diffusion of the cracks in the vault. However, it is necessary to consider also the aspects related to the materials compatibility and to the durability of the intervention. For example, among the traditional techniques adopted in the past, the use of steel-reinforced concrete plasters evidenced many drawbacks related to the introduction of excessive massive weight; moreover, the extensive use of concrete reduces the breathability of the historical masonry and steel elements evidenced, in the long term, severe problems of corrosion. Since late Nineties, the use of noncorrosive composite materials, characterized by high tensile strength to weight ratio, has gradually started as an effective alternative for the strengthening of existing masonry structures. The interaction between the reinforcement and the masonry substrate proved to be a crucial aspect, as a premature detachment at the interface can compromise the strengthening effectiveness [16]-[17]. Thus, inorganic mortars (adopted in techniques known as Textile Reinforced Mortar - TRM, Fiber Reinforced Cementitious Matrix – FRCM, and Composite Reinforced Mortar, CRM) instead of polymeric resins (adopted in the Fiber Reinforced Polymer – FRP technique), resulted preferable, due to the better fitting with rough substrates, the higher breathability, the chemical compatibility with historic masonry, the easier intervention for removal and the better performances in special environments (high or low temperatures, fire, UV-rays, water and alkaline ambient). Higher transversal loads can thus be supported by strengthened vaults, but also a higher stress level in the masonry is reached, thus other types of failure, i.e. by shear sliding or by uplift from the walls, have to be considered. In this concern, the skewback sections deserve a particular attention, for the discontinuity of the reinforcement: the connection with the abutments needs thus to be properly addressed.

The experimental results available in the literature on reinforced masonry vaults concern, mostly, monotonic or quasi-static cyclic tests based on concentrated vertical loads (at 1/3 or 1/4 of the span) [18]-[20]. As already broadly discussed by the authors [15], these simplified load patterns, although reproducing qualitatively the deformation, are not able to simulate the actual seismic behavior of non-bearing masonry vaults, as localize damage (concentrated load section) and introduce stabilizing forces at the skewback sections (vertical reactions to applied load), leading to an erroneous quantification of the reinforcement benefits. Very limited investigations were devoted also to the actual damping capacities of reinforced masonry vaults in dynamic conditions. Nevertheless, this is a crucial aspect related to the energy dissipation (through phenomena such as plastic behavior, friction, interlocking, sliding ...[21]) and to the attenuation of the amplitude of seismic motion. Corradi et al. [22] performed some dynamic identification tests on masonry barrel vaults subjected to free oscillations, comparing undamaged, damaged (formation of three hinges) and FRP retrofitted samples. The natural frequencies (1st mode) ranges around 10-13 Hz; the damping increases in the damaged vaults with respect to the undamaged ones; intermediate damping values were found for the retrofitted configurations. However, being non-destructive tests, the information on FRP retrofitted samples does not refer to near-collapse conditions. Also Giamundo et al. [23] performed some dynamic identification tests (undamaged conditions) on vaults with fixed abutments, evaluating both the unreinforced configuration and that retrofitted with TRM. The natural frequencies ranges from 13 Hz (unreinforced vault) to 19.3 Hz (retrofitted). Moreover, it emerged a damping of 2.2-3.2% for the unreinforced configuration and of 1.7%-2.8% for the retrofitted one, evidencing that, with the application of TRM, the stiffness of the unreinforced sample was restored, even slightly increased. Moreover, shaking table tests on the two configurations were carried out, evidencing the significant seismic improvement of the strengthening intervention and the progressive dissipation increase of the reinforced vault when approaching to significant damage conditions. Ramaglia et al. [24] carried out shaking table tests on vaults with moving abutments, testing both the unreinforced configuration and that strengthened through TRM at the extrados, with also the addition of reinforced masonry ribs. The evaluation in terms of damping ratio evidenced values increasing with the intensity of the signal, ranging from 2% (undamaged) to 4.5% (near failure), for unreinforced masonry, and from 1% to 10% in reinforced conditions, when approaching to significant damage. De Santis and De Felice [25] took into account a 2% damping, below the elastic threshold, and observed that reliable damping values, when damage occurs, should range from 2% to 10% (or even higher values, for very severe earthquake scenarios).

The paper deals with the strengthening of masonry vaults by means of a CRM technique based on composite Glass Fiber-Reinforced Polymer (GFRP) meshes embedded in a 30 mm thick mortar matrix applied at the vault extrados or intrados [26]. The authors already performed experimental, quasi-static cyclic tests on full-scale, isolated masonry barrel vaults lying on fixed supports, carrying their own weight and subjected to an horizontal transversal loading [15]. The tests proved that the CRM reinforcement, through the tensile strength of the GFRP wires, provides a considerable contribution in inhibiting the flexural failure in some vault sections, delaying the activation of the collapse mechanism and permitting to reach significantly high values of both resistance and displacement capacities. A simplified analytical procedure, based on Müller-Breslau equations, was also already developed to quantify the vaults lateral resistance for assigned geometry and mechanical characteristics of materials. However, the experimental and numerical studies led to concluded that, to ensure the exploitation of the whole CRM benefits, it is necessary to contrast effectively the displacement at vault skewbacks, otherwise it can fail prematurely for excessive sliding and/or uplift [27]. For this reason, a new experimental campaign was carried out, extending the study to the interaction of the vault with the masonry walls, testing a strategy to connect the CRM reinforcement with the abutments. Vaults with bricks arranged in running bond pattern and also brick *in folio* vaults (even more vulnerable to the seismic action) were tested, considering both configurations with reinforcement applied at the extrados and at the intrados. The specific vault-to-abutment connection, designed and tested, consisted in steel bars injected in the masonry walls and embedded in the mortar coating and, in case of reinforcement applied at the intrados, additional transversal GFRP connectors. The results of these new tests are presented in the paper and are analyzed also in comparison with the previous ones, allowing original considerations on the actual behavior of the connections and on the vaults dissipative capacities.

## 2    Samples characteristics

The experimental samples consisted in masonry barrel vaults made of solid bricks, having a width $w$ = 770 mm. Two groups were considered:

- the vaults of the former group (identified with prefix "v") were made of bricks arranged in a running bond pattern and had a rise-radius ratio $f/r$ = 0.75, radius $r$ = 2060 mm, span $s$ = 3830 mm and thickness $t$ = 120 mm (Figure 3a);

- the latter group (prefix "vt") concerned brick *in folio* vaults with $f/r$ = 0.5, $r$ = 2275 mm, $s$ = 3900 mm and $t$ = 55 mm (Figure 3b).

To reproduce the presence of the walls, the samples were built on masonry abutments 460 mm high and 770 mm wide. In these wallettes, the bricks were arranged in seven rows: the lower five rows, on which the vault sample lied on, had a thickness of 380 mm, the upper two rows were 250 mm thick. Before the vault construction, inclined steel shores, fixed to the concrete pavement of the laboratory, were introduced externally to avoid displacements and rotations of the masonry abutments, induced by the vault thrust. Then, the vault was built on a centering (formwork) made of polystyrene blocks. The centering was removed after at least 7 days and it was applied the reinforcement.

Six vaults samples were built, with main characteristics resumed in Table 1: the samples ID identifies the group ("v" or "vt") followed by NR, in case of unreinforced vault, or by RE, RI in case of vault strengthened at the extrados or at the intrados, respectively. Moreover, in sample vtRI+, an additional layer of mortar, 15 mm thick, was applied between the masonry and the CRM reinforcement, so to represent the cases where the removal of the existing plaster is difficult without damaging the vault and it is therefore preferable to overlay the CRM.

The solid brick units had dimensions 55x125x250 mm³ and a self-weight equal to 15.6 kN/m³ (Cov. 0.4%). Compression tests and three point bending tests [28] provided, respectively, an average compressive strength of 20.8 MPa (Cov. 4.0%,) and a flexural strength of 5.0 MPa (Cov. 4.7%). The characteristics of the masonry mortar (hydraulic lime / sand ratio = 1:3.5, in volume, self-weight 18.9 kN/m³, Cov. 1.8%), were determined through tests on prismatic [29] and cylindrical [30]-[32] samples. For each vault, characterization tests on both six mortar cylinders (three for indirect tensile tests and three for compression tests) and on three mortar prisms (three-point bending tests, with compression tests on stumps) were performed. The average compressive strength was $f_{c(c)}$ = 2.0 MPa (Cov. 23.6%), from cylinders, and $f_{c(p)}$ = 1.8 MPa (Cov. 34.6% from prisms; the tensile strengths were $f_{t(c)}$ = 0.2 MPa (Cov. 24.1%) and $f_{f(p)}$ = 0.6 MPa (Cov. 33.1%), respectively. The average Young modulus was $E_{(c)}$ = 5.9 GPa (Cov. 29.0%).

Also the characteristics of the mortar of the additional layer in sample vtRI+ (premixed lime mortar, self-weight 15.9 kN/m³, Cov. 1.0%) were determined through experimental tests (six mortar cylinders and three mortar prisms): $f_{c(c)}$ = 4.3 MPa (Cov. 5.6%), $f_{c(p)}$ = 4.5 MPa (Cov. 10.6%), $f_{t(c)}$ = 0.7 MPa (Cov. 15.1%), $f_{f(p)}$ = 2.0 MPa (Cov. 10.6%) and $E_{(c)}$ = 5.5 GPa (Cov. 6.4%).

The CRM reinforcement was applied at the vault extrados or intrados; the procedure concerned in the spray, the day before the mortar casting, of a thin layer of cement scratch coat, the drilling of holes (6 per square meter), the application of the GFRP mesh and the injection of thixotropic epoxy resin inside the holes, with the insertion of L-shaped GFRP connectors. A 30 mm thick mortar coating was then applied, so to completely embed the GFRP elements (mesh and connectors), maintaining the mesh approximately at the middle of the thickness and avoiding the formation of voids. A premixed mortar, made of natural hydraulic lime, hydrated lime and silica sand (maximum aggregate diameter 1.5mm) was used for the coating (self-weight 16.3 kN/m³, Cov. 2.3%). Characterization tests, performed on six mortar cylinders and six mortar prisms for each vault, provided average compressive strengths $f_{c(c)}$ = 8.6 MPa (Cov. 17.7%, from cylinders) and $f_{c(p)}$ = 8.2 MPa (Cov. 13.1%, from prisms). The tensile strengths were $f_{t(c)}$ = 1.1 MPa (Cov. 20.1%, splitting tests on cylinders) and $f_{f(p)}$ = 2.9 MPa (Cov. 13.8%, three-point bending tests on prisms). The average Young modulus was $E_c$ = 8.5 GPa (Cov. 11.2%). The mean mechanic characteristics of the mortars have been summarized in Table 2Table 2.

The characteristics of the GFRP elements, evaluated in accordance with [33], are summarized in Table 3. The GFRP mesh had a 66x66 mm² grid pitch. The mesh is produced by coating long glass fibers wires with a thermo-hardening resin (vinyl ester epoxy and benzoyl peroxide as catalyst) and then weaving the straight fibers wires in one direction (weft) across the wires in the perpendicular direction (warp), in which the fibers are twisted. In the vault samples, the orientation of the twisted fibers wires followed the arch direction, the parallel fibers wires the vault axis. Additional square parts of GFRP mesh (165×165mm², with 33×33mm² grid pitch) were introduced to improve the anchorage of the L-shape GFRP connectors in the mortar layer; these devices were applied over the main GFRP mesh. Actually, in sample vtRE, no GFRP connectors were used, considering that, in these thin vaults, the hole drilling could not be allowed, as might damage the valuable ceiling finishing (e.g. for the presence of *frescos*, *stuccos* or other decorations). Moreover, it is observed that, in sample vtRI+, no scratch coat was applied, but the existing mortar surface was bush-hammered before the CRM application, to improve adhesion.

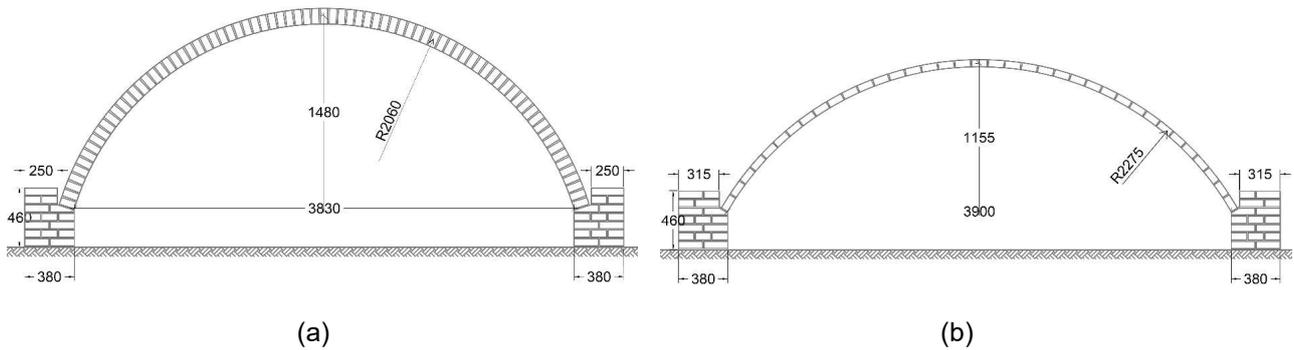

Figure 3: Main geometric characteristics of the vault samples (a) with bricks arranged in running bond pattern or (b) with bricks laid *in folio*.

Table 1. Characteristics of the tested samples: identification label *ID*, rise radius ratio *f/r*, thickness *t*, devices connecting the CRM layer with the masonry vault and with the abutments.

| ID | f/r [-] | t [mm] | Connectors CRM-vault | Connection CRM-abutments S.Steel | Connection CRM-abutments GFRP |
|---|---|---|---|---|---|
| vRE | 0.75 | 120 | 6/m² | 3φM6 | - |
| vRI | 0.75 | 120 | 6/m² | 3φM6 | 2 |
| vtNR | 0.50 | 55 | - | - | - |
| vtRE | 0.50 | 55 | - | 3φM6 | - |
| vtRI | 0.50 | 55 | 6/m² | 3φM6 | 2 |
| vtRI+ | 0.50 | 55 | 6/m² | 3φM6 | 2 |

Table 2. Average characteristics of the mortars, with respective coefficient of variation, *Cov.*: self-weight $\gamma$, compressive and tensile strengths obtained from cylinders ($f_{c(c)}$ and $f_{t(c)}$) and from prisms ($f_{c(p)}$ and $f_{t(p)}$) and Young modulus $E_{(c)}$.

| Mortar | γ [kN/m³] | Cov. [%] | Cylinders $f_{c(c)}$ [MPa] | Cov. [%] | $f_{t(c)}$ [MPa] | Cov. [%] | $E_{(c)}$ [GPa] | Cov. [%] | Prisms $f_{c(p)}$ [MPa] | Cov. [%] | $f_{t(p)}$ [MPa] | Cov. [%] |
|---|---|---|---|---|---|---|---|---|---|---|---|---|
| Masonry mortar | 18.9 | 1.8 | 2.0 | 23.6 | 0.2 | 24.1 | 5.9 | 29.0 | 1.8 | 34.6 | 0.6 | 33.1 |
| Mortar for additional layer | 15.9 | 1.0 | 4.3 | 5.6 | 0.7 | 15.1 | 5.5 | 6.4 | 4.5 | 10.6 | 2.0 | 10.6 |
| Mortar for the coating | 16.3 | 2.3 | 8.6 | 17.7 | 1.1 | 20.1 | 8.5 | 11.2 | 8.2 | 13.1 | 2.9 | 13.8 |

Table 3. Characteristics of the GFRP elements: dry fiber and equivalent cross sections, $A_{fib}$ and $A_{eq}$, tensile resistance, $T_w$, ultimate elongation, $\varepsilon_{ult}$, and axial stiffness, $EA_w$, with respective coefficient of variation, *Cov.*

| GFRP element | $A_{fib}$ | $A_{eq}$ | $T_w$ | Cov. | $\varepsilon_{ult}$ | Cov. | $EA_w$ |
|---|---|---|---|---|---|---|---|

|  | [mm²] | [mm²] | [kN] | [%] | - | [%] | [kN] |
|---|---|---|---|---|---|---|---|
| Straight fibers wire | 3.8 | 9.41 | 5.62 | 4.8 | 0.0193 | 5.7 | 291 |
| Twisted fibers wire | 3.8 | 7.29 | 4.49 | 6.7 | 0.0169 | 9.5 | 266 |
| Connectors | 57.6 | 7x10 | 20.96 | 10.9 | - | - | - |

## 3 Focus on the CRM connection with the abutments

The specific requirements for the design of the connection of the reinforced vaults with the abutments are individuated in Figure 4: to guarantee an adequate support, sliding in both the transversal (Figure 4a) and in the longitudinal direction (Figure 4b) has to be prevented, as well as the uplift (Figure 4c). However, for the intervention optimization, some rotation at the skewback section should be allowed: this permits the vault to deform significantly and to dissipate energy through the widespread cracking of the CRM layer.

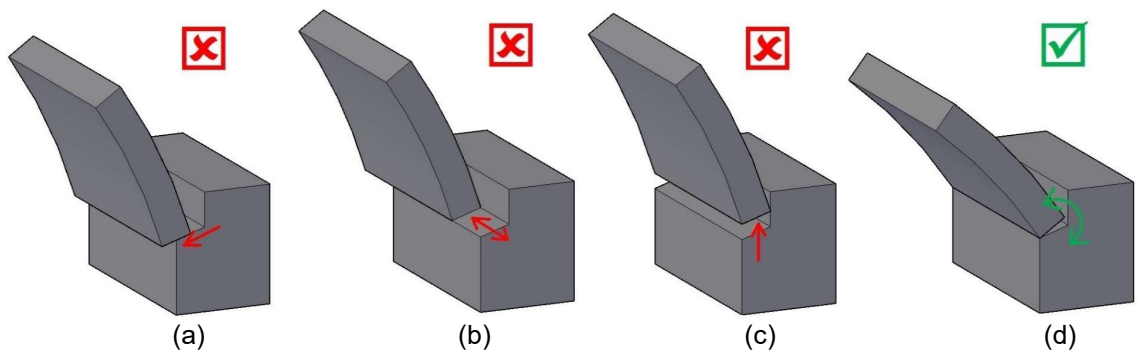

Figure 4: Requirements for the design of the connection of the reinforced vaults with the abutments: contrast the sliding at the skewbacks in the (a) transversal and (b) longitudinal direction and (c) the uplift; (d) permit some rotation.

Actually, the outward transversal slip at the vault skewbacks is avoided by the upper rows of the masonry wallettes, representing the walls (the applied constraints are described in §4). Differently, to prevent the inward sliding, tensile resistant elements connecting the CRM layer with the abutments were introduced in the configuration reinforced at the extrados. A higher resistance against the inward shear action was thus indirectly provided to the cracked skewbacks by increasing the stress on the compressive side (the shear resistance is mainly due to friction). Steel bars were thus introduced; in particular, due to the thin mortar cover, stainless steel was preferred, so to avoid durability drawbacks related to corrosion. M6 threaded bars, AISI 316, were considered (net cross section 20.1 mm$^2$). The tensile characterization tests performed on three samples [34] (Figure 5.a-b) provided an approximately elastic-perfectly plastic behavior (dotted lines curves in Figure 6), with mean axial stiffness $EA_b$ = 3417 kN (Cov. 13%) and resistance $T_b$ = 16.4 kN (Cov. 1.0%). Moreover, to exploit the ductility of the steel bars, allowing some rotation at the vaults skewbacks (Figure 4d), their resistance has to be lower than that of the GFRP wires (which, on the contrary, have a brittle behavior). Thus, considering a GFRP grid pitch of 66 mm wuth a tensile strength of 4.49 kN for the single GFRP wire (Table 3), three M6 bars were introduced at both skewbacks (global tensile strength of the bars 3 * 16.4 kN = 49 kN < 52.4 kN = 4.49 kN * 770 mm / 66 mm). Obviously, to allow the bars to yield, it is necessary to guarantee an adequate bond length. Pull-out tests were performed to evaluate the behavior of M6 bars embedded in a 30 mm thick CRM layer applied on a masonry wall (Figure 5c-e). The load, $T$, was applied by a hand pump through a hollow hydraulic cylinder (200 kN capacity, 50 mm stroke) equipped with a pressure transducer. The displacement between the bar and the mortar coating edge, $\Delta l$, was monitored by two displacement transducers (stroke 50 mm); the initial measurement base length (the free length from the jacketing) was 112 mm. The load-displacement graphs of the pull-out tests referred to a bond length equal to 270 mm (45 times the bar nominal diameter) are plotted in Figure 6, in comparison with those obtained from the tensile tests, in which the initial gauge length was the same. The two sets of curves overlapped in the first part of the tests, but then gradually diverged from about 5-6kN, likely due to the progressive detachment of the bar from the mortar. A sudden load decrease then occurred as the bar, once completely detached form the mortar, was extracted. The pull-out tests attained to a mean load value of 15.2 kN (Cov. 3.2%), which resulted very close to yielding strength of the bar (Figure 6). It can thus be deduced that the provided bond length is approximately the limit value to be ensured. It was therefore decided, for the vault samples, to set the bond length to 320 mm (50 times the bar nominal diameter) and to apply an end nut to the bar so to provide further residual contrast to its extraction from the mortar.

The detail of the skewback configuration for the vaults reinforced at the extrados is reported in Figure 7, for both the 120 mm and the 55 mm thick samples. Each bar was injected in the masonry abutment with tixotropic epoxy resin (300 mm bond length) and embedded in the mortar coating (320 mm bond length). To be noted that the GFRP mesh is placed above the bars, to better counteract the possible mortar splitting (the bonding of the steel bars is guaranteed by the indentation of the threads in the mortar, which causes tensile stresses in the mortar in the direction perpendicular to the bar axis and may induce the mortar crack in the direction of the bar). The bars are able to contrast also longitudinal sliding and uplift.

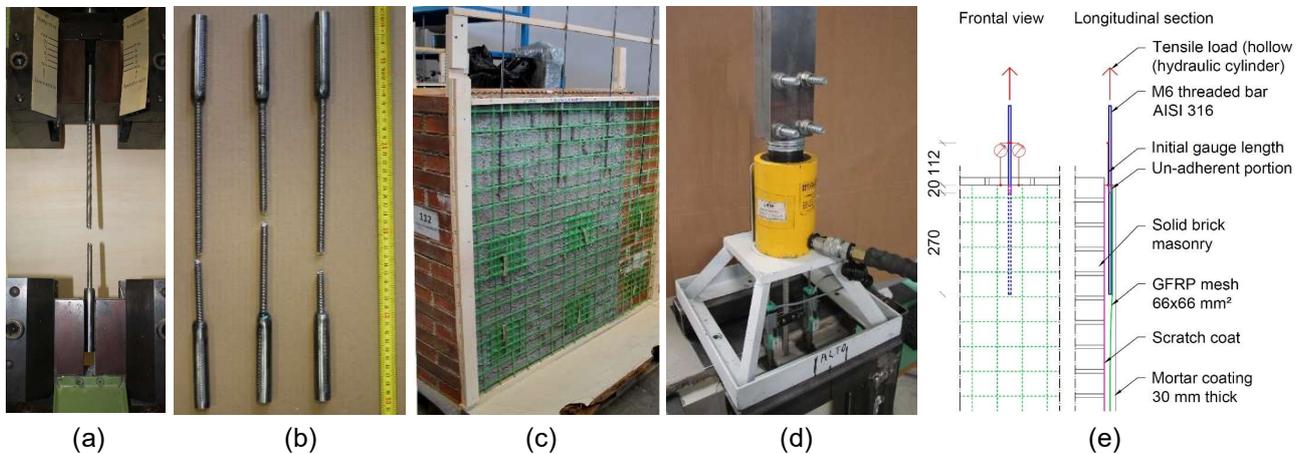

(a)          (b)          (c)          (d)          (e)

Figure 5. Characterization of M6 stainless steel threated bars: (a, b) tensile tests and (c, d, e) pull-out tests.

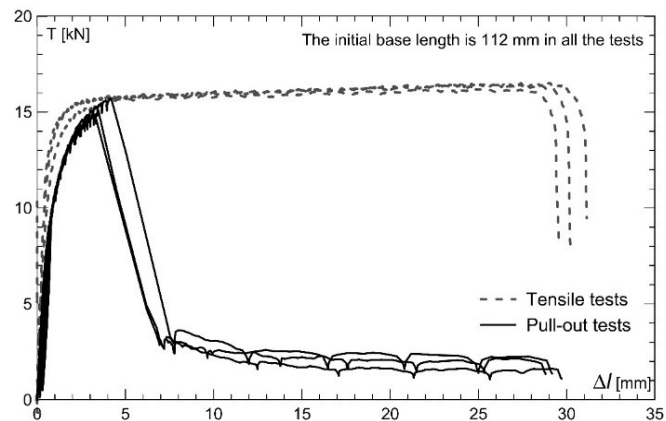

Figure 6. Load-displacement graphs of M6 stainless steel threated bars subjected to tensile tests (dotted lines) and to pull-out tests from CRM layer, with bond length 270 mm (solid lines).

Also the configurations with the CRM reinforcement applied at the intrados were provided with 3 M6 steel bars at each skewback, preventing longitudinal sliding and uplift (Figure 4b-c) and permitting some rotation (Figure 4d). Moreover, it was necessary to introduce further elements to prevent the inward sliding (Figure 4a), as it is not adequate to rely only on the dowel effect of these bars. In particular, the transversal L-shape GFRP connectors closest to skewback sections were deeply injected in the masonry abutments (Figure 8); two GFRP connectors at each skewback section thus resulted. Two pull-out tests were carried out to evaluate the behavior of these elements: the samples consisted in a 380x380x30 mm³ mortar slab with the GFRP mesh and a single L-shaped GFRP connector embedded; also the 33 x 33 mm² piece of GFRP mesh was utilized (Figure 9a). The sample was placed over a rigid flat support and the quasi-static extraction force was applied to the free connector edge through two hydraulic jacks (150 kN, 101 mm stroke) connected in parallel and activated by a hand pump (Figure 9b). During the tests, as the load increased, a mortar splinter cracked close to the connector embedded edge (Figure 9c). Then the load was almost maintained since the connector was gradually extracted from the slab. The average failure tensile load was $T_c$ = 5.6 kN (Cov. 12.8%). So, at the skewback section of the vaults reinforced at the intrados, it was awaited an inward shear resistance of about 11.2 kN, due to the presence of 2 GFRP connectors (= 2·5.6 kN).

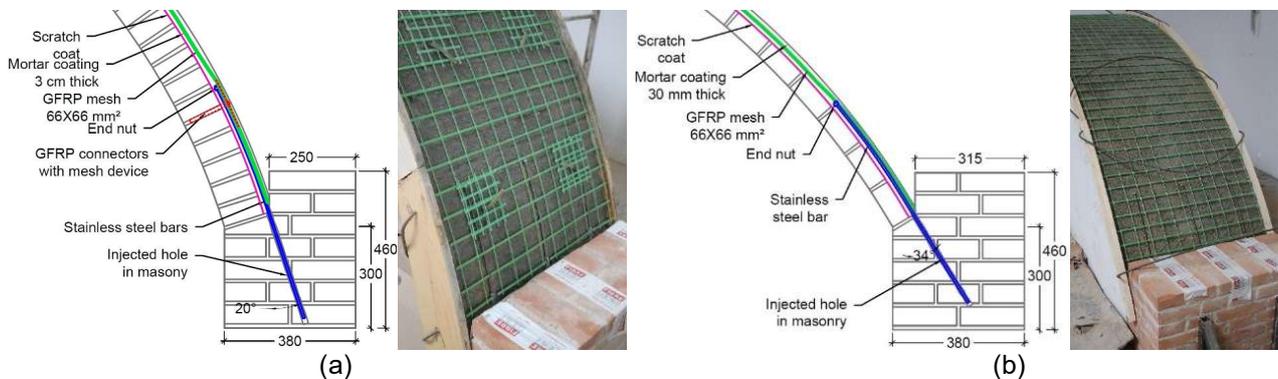

Figure 7. Abutments configuration, in case of CRM reinforcement applied at the extrados: samples (a) vRE and (b) vtRE.

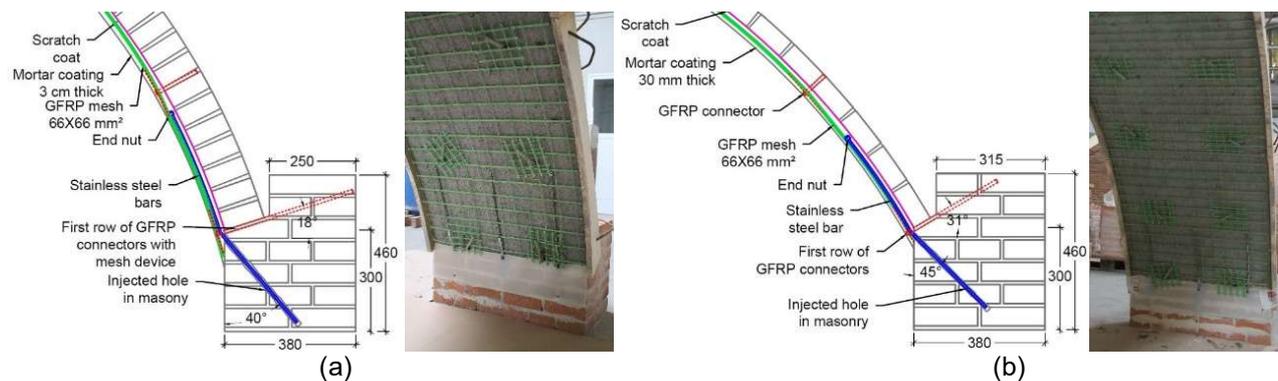

Figure 8. Abutments configuration, in case of CRM reinforcement applied at the intrados: samples (a) vRI and (b) vtRI.

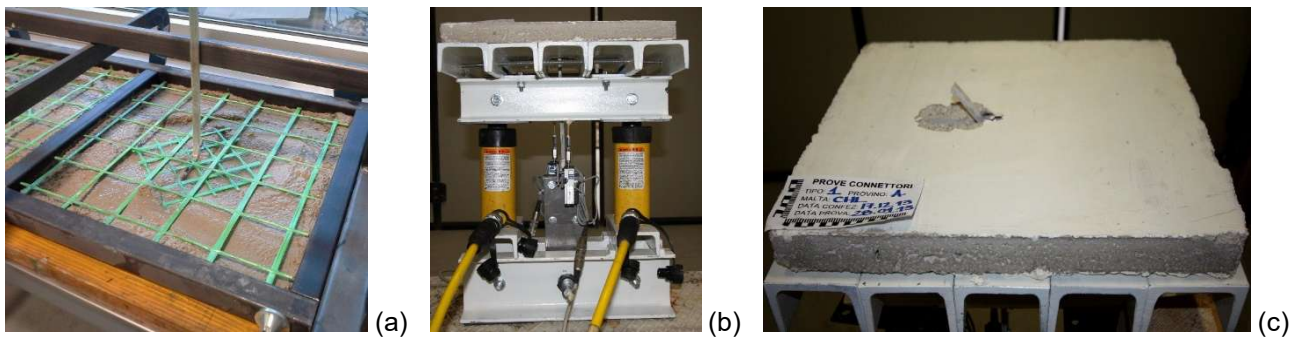

(a)                                              (b)                                              (c)

Figure 9. Extraction tests of the GFRP connectors from the CRM layer: (a) sample construction, (b) test setup and (c) example of the configuration at the end of the test

## 4    Experimental setup

The test setup (Figure 10a) was aimed to reproduce a transversal load pattern proportional to the vault mass so to represent the effect of the horizontal inertial forces induced by an earthquake on masonry vaults lying on fixed abutments and supporting their own weight only. Quasi-static cyclic tests were carried out: forces of identical magnitude and direction were applied at eight different points along the vault length, equally spaced, through pinned frame systems installed at the two halves of the vault and supported by sliding guides. The horizontal load was modulated through two hydraulic actuators connected in parallel, applied at the opposite sides of the contrast structure and activated by a hand pump. Flow control valves permitted loading, unloading and the inversion of the loading direction. The steel apparatus was designed adequately stiff, so to assume negligible its deformation in respect to that of the sample; loading cells and potentiometer transducers allowed the real-time monitoring of the load/displacement story. The horizontal displacements of the masonry abutments were constrained through transversal steel beams bolted to the main apparatus (Figure 10b). Moreover, a vertical pre-compression was applied, reproducing the axial load transmitted by the upper wall and avoiding the abutment failure for horizontal shear. Thin layers of gypsum were applied along the masonry surface and the front faces, so to facilitate to view the occurrence of the cracks during testing.

More details were reported in [15], as the test apparatus and method were the same as those adopted in previous tests. However, it is observed that, due to the thinness of the brick *in folio* vaults, the vertical forces transmitted by the loading arms, induced by their weight, can modify appreciably the vault stress state. Thus, in this new testing campaign, the loading arms were hanged on the upper beams of the contrast frame through steel cords equipped with a preloaded tension spring at the upper end (Figure 10b). The initial elongation of each spring was calibrated so to support the vertical load transmitted to the vault. Similarly, also the pinned frame systems for the load distribution were hanged to avoid friction on the sliding guides. In the early stages of the tests, the variation in the springs elongation is negligible, since the vault deformation is limited, as well as the cords inclination (the cords length was at least 950 mm).

The scope of the tests was to comprehend how the CRM reinforcement modifies the resistance mechanisms of the vaults against lateral loads. Clearly, the simplified laboratory asset was not able to account for the dynamic variability of the seismic action but, realistically, the resistance mechanism of the CRM reinforced masonry vault is analogous: basically, the high tensile strength of CRM provides bending resistance to the vaults sections. Moreover, similar resistance mechanisms reasonably activate also in vaults lying on moving abutments.

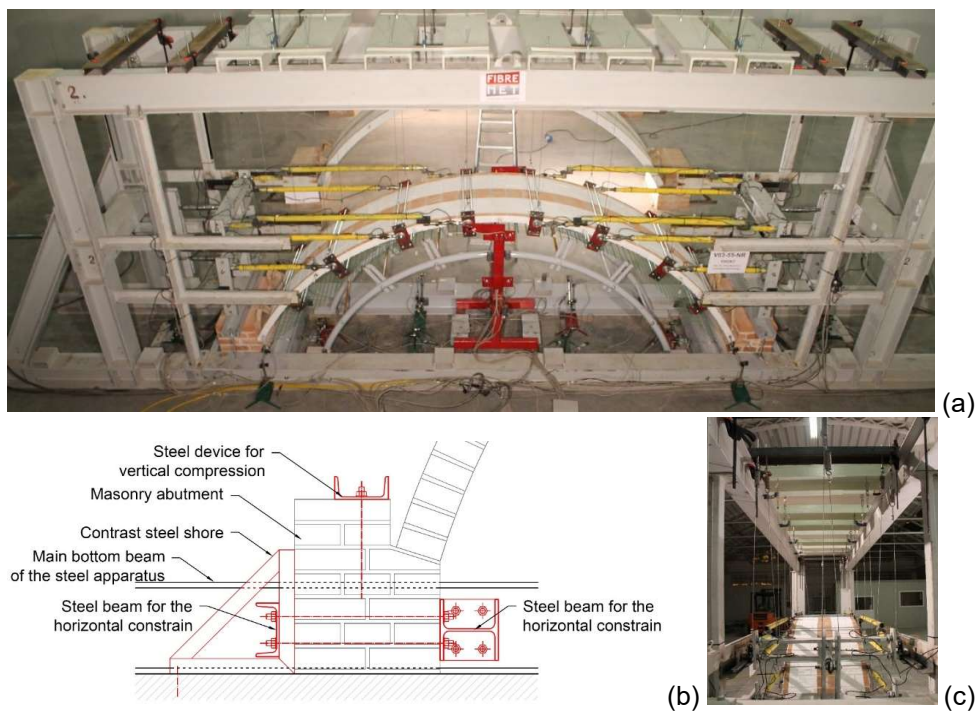

Figure 10: Test setup for the vaults: (a) global view and details (b) of the abutment arrangement and (c) of the hanging system

## 5    Tests results

The performances of the CRM reinforced vaults are reported in Figure 11a-d, as capacity curves expressed in terms of global horizontal load, $F_h$, against the mean horizontal displacement monitored at keystone $\delta_h$.

Moreover, the crack patterns of the vaults are schematically represented in Figure 12a-e, distinguishing cracks formed for the two opposite loading directions (red and blue colours in agreement with the load direction). In general, in the test on reinforced samples, the two skewback sections and two opposite haunch sections cracked in alternate position intrados-extrados. However, after the attainment of this "first significant damage" configuration, the GFRP mesh and the steel bars connection with the abutments avoided the opening of structural hinges in the sections where the reinforcement was located on the side subjected to tension. Subsequently, in the "post-damage" stage, a gradual stiffness degradation was observed, as cracks widespread in the mortar coating at the two haunches, alternately (Figure 13b and Figure 14b). At the opposite haunch, several cracks formed also on the masonry side, mainly in the bed joints (Figure 13c and Figure 14c).

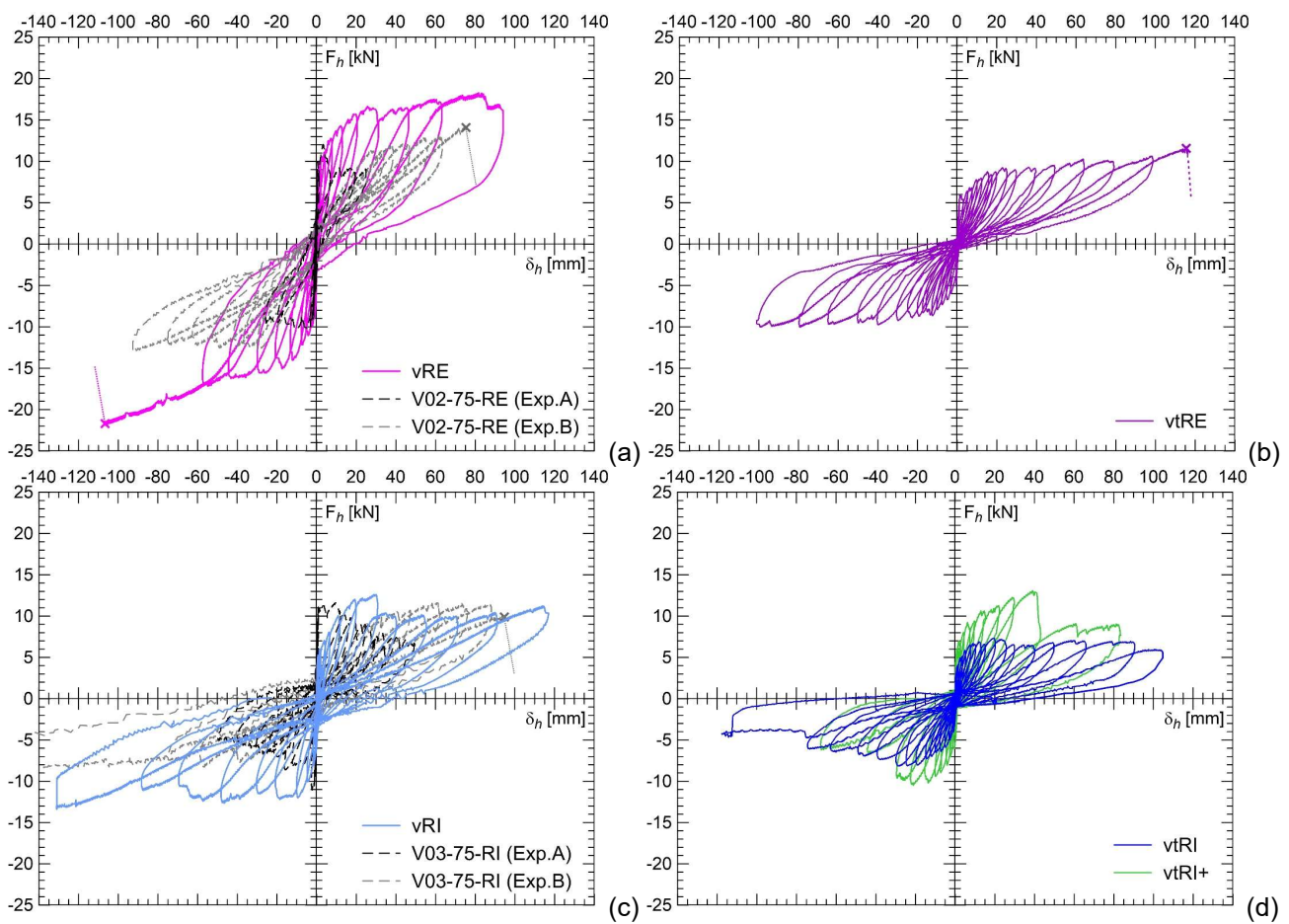

Figure 11. Capacity curves of samples (a) vRE, (b) vtRE, (c) vRI, (d) vtRI and vtRI+. The capacity curves of the previous experimental tests on samples V02-75-RE and V03-75-RI are also reported in (a) and (c), for comparison.

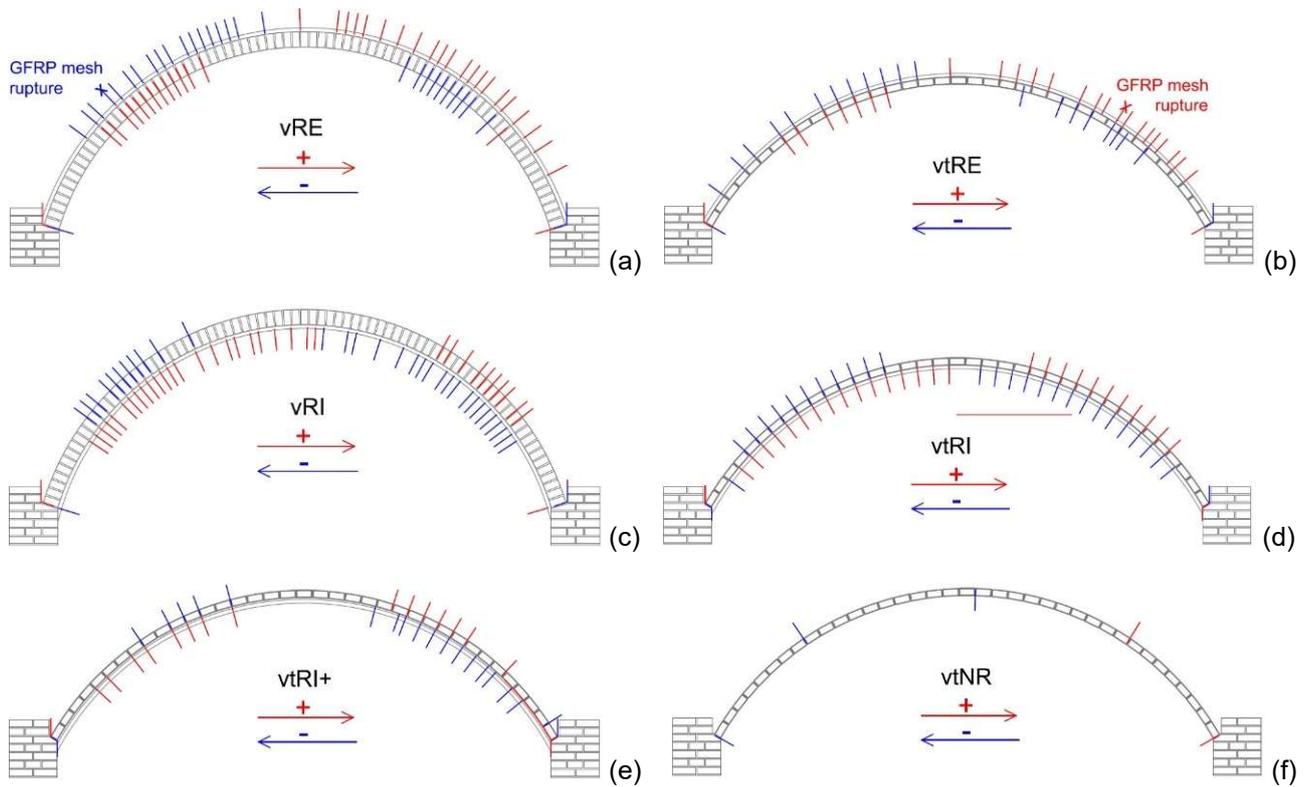

Figure 12. Crack pattern of samples (a) vRE, (b) vtRE, (c) vRI, (d) vtRI, (e) vtRI+ and (f) vtNR

In the vaults reinforced at the extrados, vRE and vtRE (Figure 13), at the increasing of the deformation, the load kept increasing. The collapse was due to the rupture of the GFRP wires in the vault in correspondence of a haunch section (at the extrados). No evident slips at the skewbacks were noted during the tests (Figure 13e), but conic splinters in correspondence of some steel bars, close to the end nuts, were observed, as likely caused by their partial extraction from the coating (Figure 13d). Sample vRE reached a peak value $F_{h,max}$ of 22.9 kN, in correspondence of a displacement $\delta_{h,u}$ = 108.1 mm; sample vtRE reached a maximum load $F_{h,max}$ = 11.8 kN at $\delta_{h,u}$ = 115.5 mm.

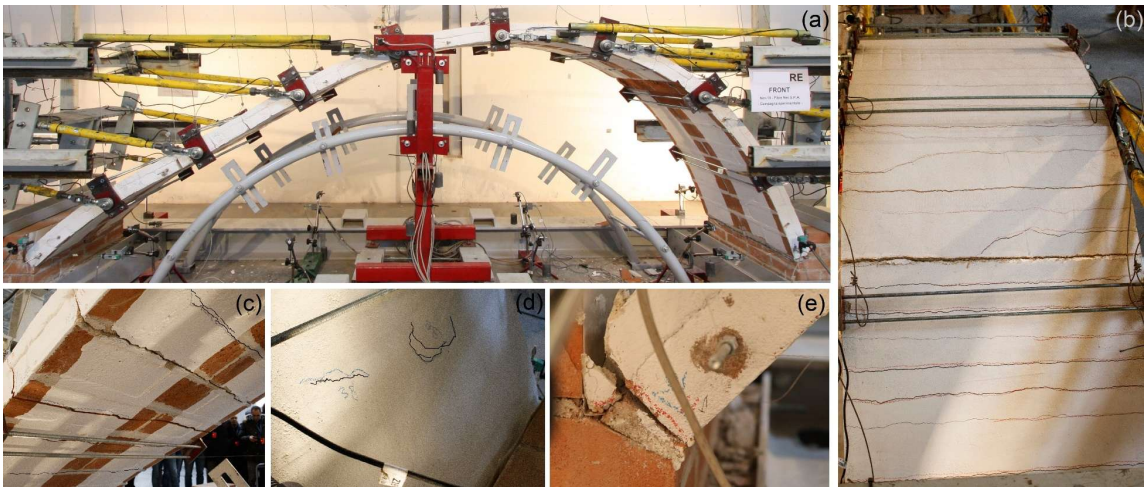

Figure 13. Vault samples reinforced at the extrados, at the incipient collapse (example of vtRE): (a) global view and typical crack pattern details concerning (b) diffuse cracking of the mortar coating at the extrados, right haunch, (c) cracking of the masonry at the intrados, left haunch, (d) local cracking of the mortar in correspondence of the connection steel bars and (e) rotation in correspondence of the left skewbacks, with no evident inward slip.

Differently, in the samples reinforced at the intrados (Figure 14), vRI, vtRI, and vRI+, some inward slip occurred (Figure 14e), as evidenced also by the partial pull-off (mortar splintering) of the first row of GFRP connectors in sample vRI (Figure 14d). Moreover, a conic mortar splinter in correspondence of some steel bars, close to the end nut, detached. These damages in correspondence of both the GFRP connectors and the steel bars did not allow further load increases. In addition, in vtRI+, a detachment surface formed between the pre-existing mortar layer and the masonry substrate, starting from the right skewback and reaching the haunch section (Figure 15). The peak load of sample vRI ($F_{h,max}$ = 13.4 kN), was reached at the ultimate cycle ($\delta_{h,u}$ = 130.9 mm); then the test was interrupted due to the considerable displacement achieved. However, at these displacement values, some signals of the GFRP wires incipient breakage were clearly perceived. Differently, sample vtRI attained to +7.3 kN at $\delta_h$ = +19.6 mm and to -8.2 kN at $\delta_h$ = -27.9 mm; then the load was approximately maintained in the positive direction, till about +105.3 mm, while gradually decreased in the negative one, due to the accidental, premature damaging of one GFRP connector, which emphasized the inward sliding. Sample vtRI+ attained to -10.5 kN (at $\delta_h$ = -21.0 mm) and to +13.1 kN (at $\delta_h$ = +39.0 mm), when the detachment of the pre-existing plaster layer formed and the load suffered an abrupt partial drop down.

In the unreinforced vault, vtNR, the cracks occurred in few sections (Figure 12f). At first, in the negative loading direction, three cracks opened in the left vault side: one at extrados, at a haunch section, and two at the intrados, at the skewback and near the keystone. Then, when the loading direction was inversed, one crack at the extrados (at right haunch) and one at intrados (at the right skewback) formed. The presence of at least four hinges, in alternate position extrados-intrados, was thus achieved and the vault failure mechanism activated. The maximum load, $F_{h,max}$ = 0.55 kN, was attained in correspondence of an horizontal displacement at the crown section $\delta_{h,u}$ < 1 mm. So, due to the very low performances in terms of both resistance and displacements, the capacity curve of the unreinforced masonry sample was omitted.

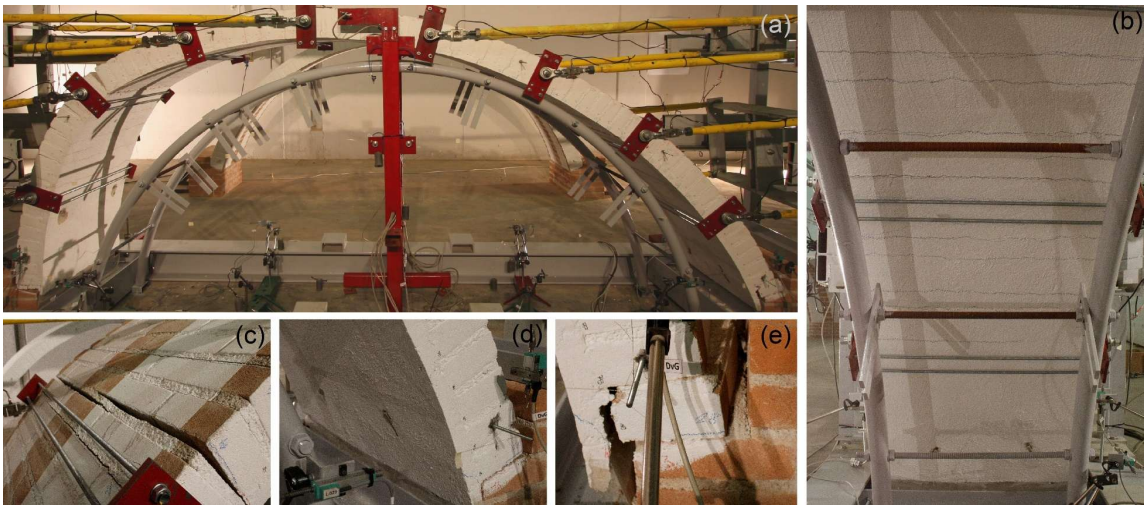

Figure 14. Vault samples reinforced at the intrados, at the incipient collapse (example of vRI): (a) global view and typical crack pattern details concerning (b) diffuse cracking of the mortar coating at the intrados, right haunch, (c) cracking of the masonry at the extrados, left haunch, (d) local mortar pull-off in correspondence of the first row of GFRP connectors and (e) rotation in correspondence of the right skewbacks, evidencing some inward slip.

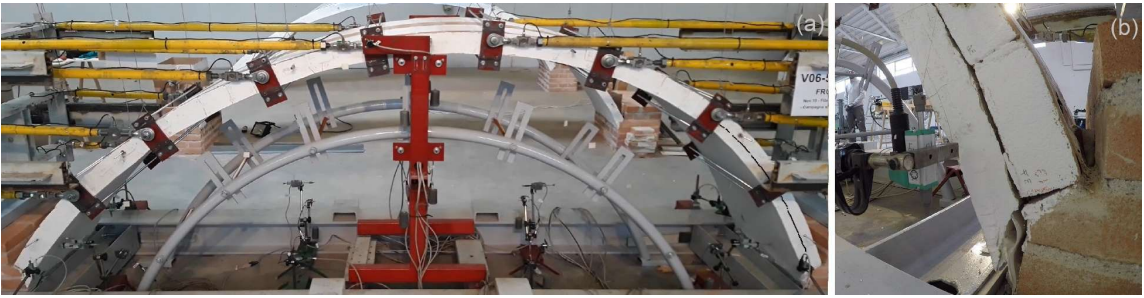

Figure 15. Configuration of sample vtRI+ at the incipient collapse, evidencing the partial detachment of the pre-existing plaster layer from the masonry, on the right side: (a) global view and (b) detail at the right skewback.

To allow comparisons, also the results of some previous tests, detailed and described in [15], are reported in Figure 11a, c: in sample V02-75-RE the CRM reinforcement was applied at the extrados, in V03-75-RI at the intrados. The samples had the same geometric and mechanic characteristics of type "v" vaults, but were built of fixed steel supports and the CRM reinforcement was interrupted in correspondence of the skewback sections (Figure 16a, c). In the first part of the tests on strengthened samples (dotted black curves "Exp.A" in Figure 11a, c), just after the bending cracking at the vault skewbacks, some slide at the skewbacks occurred, inducing a significant resistance decrease. In fact, the high lateral load produced a consistent shear action at the vault skewbacks, but reduced significantly the compression at one skewback (and, thus, its shear resistance, mainly due to friction).

In particular, the outward slip of the vault skewbacks was observed in the sample V02-75-RE; so, external steel profiles were introduced to effectively oppose it, reproducing the presence of a masonry wall (Figure 16b). This allowed, in the second phase of the test (dotted grey curve "Exp.B" in Figure 11a), the load recovery and the exploitation of the CRM resistance, till the GFRP mesh failure at an haunch section (maximum horizontal load $F_{h,max}$ = 14.1 kN, at $\delta_{h,u}$ = 76.0 mm). Moreover, it is observed that the inward slip at skewback sections did not occur, as it was contrasted by the presence of the end portion of the CRM layer, in the fixed bricks row (Figure 16a).

In sample V03-75-RI, besides the installation of the external steel profiles, polyester belts were introduced to contrast also the inward slip (Figure 16d). Actually, a lower belt tightening at the right skewback determined less restraint effectiveness in the negative loading direction, resulting in an unsymmetrical capacity curve (dotted grey curve "Exp.B" in Figure 11c). The GFRP mesh failure at an haunch section was attained at $F_{h,max}$ = 11.6 kN ($\delta_{h,u}$ = 92.9 mm).

One additional specimen was tested unreinforced (V01-75-NR) and collapsed at $F_{h,max}$ = 3.8 kN ($\delta_{h,u}$ = 0.84 mm) for the loss of equilibrium in consequence of the sudden formation of four structural hinges, induced by the masonry flexural failure for tensile cracking, in an alternate position intrados-extrados.

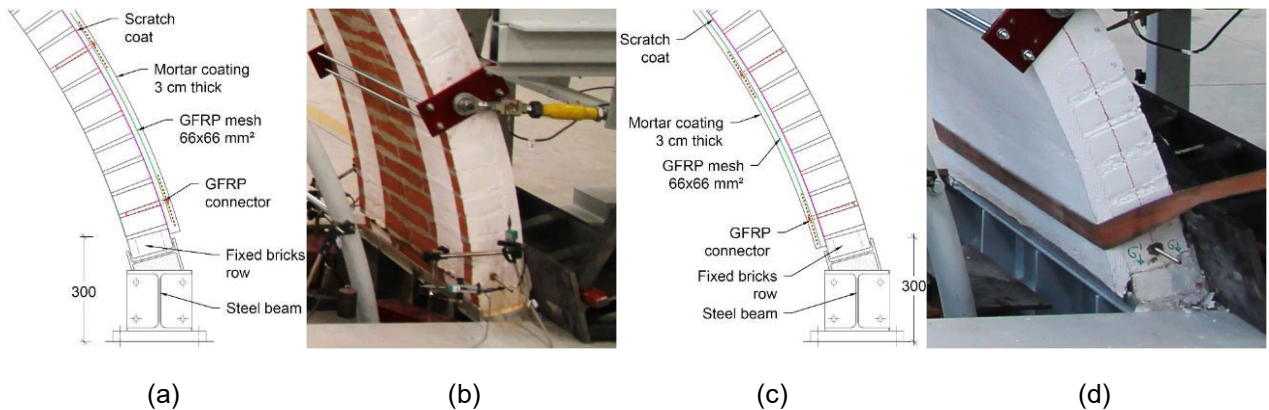

|   (a)   |   (b)   |   (c)   |   (d)   |

Figure 16. Previous tests on masonry vaults: detail of the skewback configuration (a) in V02-75-RE, with (b) introduction of the external steel profile to contrast outward sliding, and (c) in V03-75-RI, with (d) introduction of the polyester belt to contrast inward sliding

# 6    Analysis and discussion of the results

The results are summarized in Figure 17 in terms of peak load, $F_{h,max}$, and ultimate horizontal displacement $\delta_{h,u}$ at the crown section (the dotted line columns and the samples ID in round brackets are referred to the three specimens tested in the previous experimental campaign).

In general, it is worth to note that the 30 mm mortar layer induced an increase in the initial vault horizontal stiffness (about 50% in single wythe vaults and about 100% in the thinner ones). An increase in the stiffness (and, therefore, in the main eigenfrequency) of a structure may imply an increase in the seismic action, but due to the high stiffness of the unstrengthened vault [22]-[23], the spectral amplification affecting the reinforced vault is at most the same. In cracked conditions, the stiffness and the capacity of unreinforced vaults significantly reduced up to the formation of the collapse mechanism. Differently, in CRM reinforced vaults, due to the presence of the GFRP mesh that opposed to crack opening, the stiffness reduced gradually, the capacity was maintained/increased and very high ultimate displacement values were reached.

The resistances of samples vRE and vRI resulted, respectively, 6.1 and 3.6 times that of the unreinforced sample V01-75-NR; moreover, ultimate displacements more than 100 times greater were attained. As expected, sample vRE provided higher performances in respect to V02-75-RE (Figure 11a): this was due to the additional rotational constraint provided by the connection with the abutments. In fact, the steel bars, acting in tension, opposed to the free rotation when the skewback section cracked. The behavior of sample vRI resulted similar to that of V03-75-RI supplied with an adequately tightened belt (positive loading direction) - Figure 11c: the steel bars provided some rotational constraint but the GFRP shear connector allowed some inward sliding.

Reasonably, lower resistance values were generally achieved by the brick *in folio* vaults with the same reinforcement configurations, due to the thinner cross section and higher slenderness, despite the reduced *f/r* ratio. However, the reinforcement effectiveness, evaluated in terms of ratio between the peak load of the reinforced and the unreinforced vaults, was significantly greater: 14.9 for vtRI and 21.5 for vtRE.

The comparison between vtRI+ and vtNR is not proper (an unreinforced masonry vault with an additional mortar layer should be considered). On the contrary, samples vtRI and vtRI+ can be compared: the presence of the additional layer of mortar in vtRI+, led to a resistance 60% higher than vtRI, due to the thicker cross section; displacement capacity resulted lower because of the detachment of the mortar from the masonry substrate (-63%).

Considering the results of the characterization tests carried out on the connections with the abutments (section 3), greater deformations actually occurred at the reinforced vaults skewbacks. In particular, the partial sliding of the steel bars from the mortar allowed some additional rotation at the skewbacks. Moreover, some inward slip at the skewbacks occurred in the configurations with CRM at the intrados, as evidenced also by the local damaging of the mortar around the GFRP connectors. Likely, the complexity in the application of the CRM on the vaults and the presence of a cyclic action influenced the performances of the connection devices, reducing their effectiveness in respect to that obtained from the characterization tests. However, in respect to the unreinforced configuration, very high performance improvements always resulted from the CRM reinforced configurations and the loss of the support at the skewbacks was always prevented. The application of CRM at the extrados resulted more effective than that at the intrados.

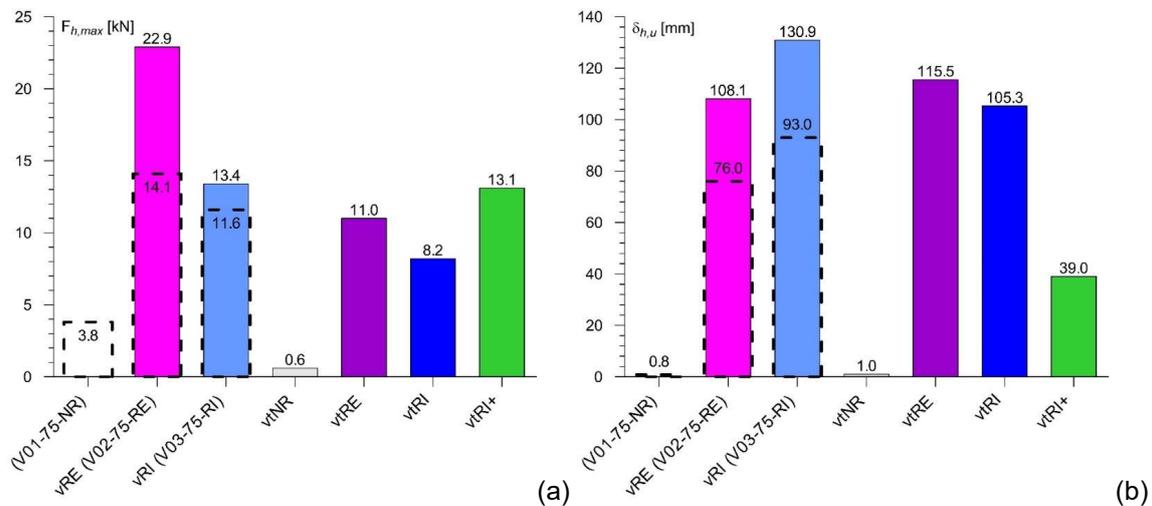

Figure 17. Experimental results in terms of (a) lateral resistance $F_{h,max}$ and (b) ultimate displacement $\delta_{h,u}$. Dotted line columns refer to the results of the three previous tests (samples ID in round brackets) [15].

The quasi static-cyclic tests provided also useful information on the dissipative performances of the reinforced vaults, providing indications on the damping to consider for the evaluation of the design response spectrum [35] for the CRM reinforced vaults. In particular, an evaluation of the energy dissipation capacity of the samples was carried out by assessing the input energy $\Delta E_{inp}$ and the dissipated hysteretic energy $\Delta E_{hys}$ of each cycle. In particular, $\Delta E_{inp}$ was calculated as the work made by the actuators to reach the target displacement amplitude $\delta_{h,T}$ (both in the positive and negative loading direction) and $\Delta E_{hys}$ as the area included in the hysteretic loop – Figure 18a. The trends of $\Delta E_{inp}$ and $\Delta E_{hys}$ are shown in Figure 18b. In particular, the attention focused on the behavior of the reinforced vaults after the attainment of the target displacement $\delta_{h,T}$ associated to "first significant damage" configuration (cracking of both the skewback sections and the haunch sections). It was found that, in the "post-damage" phase, the ratios between $\Delta E_{hys}$ and $\Delta E_{inp}$ remained almost constant in all cycle up to vault collapse, as resumed in Table 4. Similarly for the coefficient of hysteretic damping $\beta_0$, which was calculated for each cycle, through the following equation [21], [35]:

$$\beta_0 = \frac{1}{2\pi} \frac{\Delta E_{hys}}{\Delta E_{S0}},$$ (1)

where the strain energy, $\Delta E_{S0}$, was individuated by the area of the right triangles having a vertex in the axis origin and a vertex in correspondence of the target point (evidenced by the dashed lines in Figure 18a).

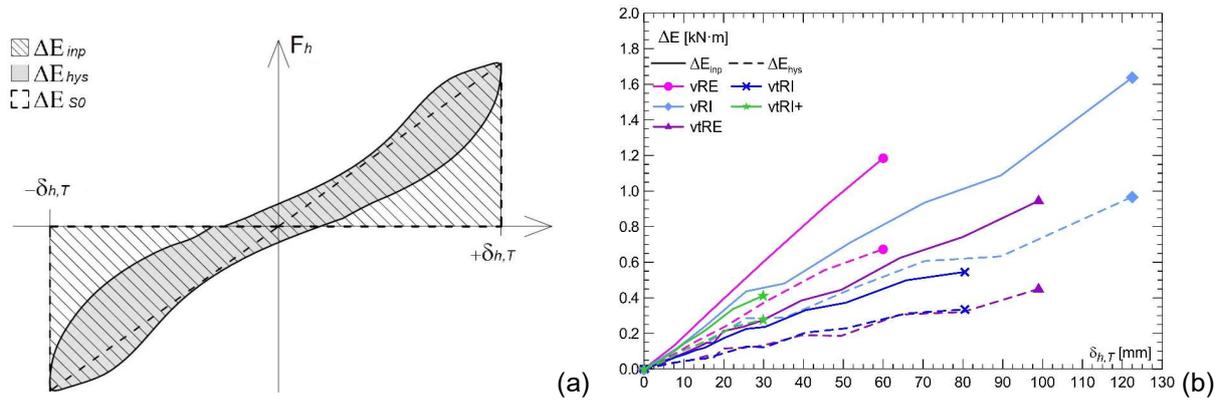

Figure 18. Dissipative performances of the samples: (a) individuation of the input energy $\Delta E_{inp}$, of the hysteretic energy $\Delta E_{hyp}$ and of the strain energy $\Delta E_{S0}$ in a cycle and (b) trends of $\Delta E_{inp}$ and $\Delta E_{hyp}$ at the increasing of the target displacement.

The ratios $\Delta E_{hys}/\Delta E_{inp}$, as well as $\beta_0$, resulted very similar for both vRE and vRI, with average values of 0.61 and 12.4%, respectively. Differently, some variations were observed in the brick *in folio* vaults: lower values emerged for the sample reinforced at the extrados, vtRE (0.48, 8.1%), in respect to vtRI (0.56, 10.4%). Moreover, the presence of the pre-existing mortar layer, in vtRI+, determined higher dissipation performances (0.72, 16.3%), compensating the lower displacement capacity.

According to Capacity Spectrum Method based on the Equivalent Viscous Damping reported in FEMA440 [35], the dissipative capacity of a structure, to be considered for determining the design response spectrum, can be evaluated considering the equivalent damping $\beta$ - eq. (5):

$$\beta = \beta_{el} + \beta_0$$ (5)

It accounts for the inherent damping $\beta_{el}$ (due to local, inelastic phenomena at the micro-scale level, which occur also below the global elastic threshold) and the hysteretic damping $\beta_0$ - eq.(1). As emerged from the dynamic investigations discussed in section 1, $\beta_{el}$ should be around 1-2%. Thus, the equivalent damping of the investigated CRM vaults in the post-damage phase reasonably ranges around a mean value of 13%.

Table 4. Mean values of the dissipated energy ratio $\Delta E_{hys}/\Delta E_{inp}$ and of the coefficient of hysteretic damping, $\beta_0$, with respective standard deviations, *St.dev.*, in the post-damage phase (after the attainment of the target displacement $\delta^*_{h,T}$)

| Sample | $\delta^*_{h,T}$ [mm] | $\Delta E_{hys}/\Delta E_{inp}$ | *St. dev.* | $\beta_0$ [-] | *St. dev.* |
|---|---|---|---|---|---|
| vRE | ±7.5 | | 0.02 | 12.5% | 1.3% |
| vRI | ±8.0 | 0.61 | 0.03 | 12.3% | 1.3% |
| vtRE | ±5.5 | 0.48 | 0.04 | 8.1% | 1.4% |
| vtRI | ±5.3 | 0.56 | 0.05 | 10.4% | 1.2% |
| vtRI+ | ±7.5 | 0.72 | 0.03 | 16.3% | 1.7% |

## 7 Conclusions

An elevate seismic vulnerability affects secondary vaulted elements in historic masonry buildings, evidencing the necessity to investigate on effective and compatible strengthening solutions. The original results of a recent experimental campaign on CRM strengthened vaults, carrying their own weight and subjected to quasi static transversal cyclic load are reported. The reinforcement technique was based on the application, at the vault extrados or intrados, of a 30 mm thick layer of mortar with GFRP meshes embedded. A particular attention was devoted to the role of the connections between the vault and the abutments, as in previous studies the authors had evidenced the importance of this detail to exploit the CRM reinforcement benefits. In fact, an ineffective connection leads to displacements at the vault skewbacks which can also induce its premature collapse for loss of the support.

A solution to design and implement this connection in practice is proposed in the paper: stainless steel bars are embedded in the mortar of the coating and injected at the masonry abutments; moreover, in case of CRM reinforcement applied at the intrados, GFRP "L-shape" connectors, injected at the masonry abutments, were introduced at the skewbacks to provide additional contrast against inward slip. Actually, the outward slip is avoided by the masonry walls, whose presence was reproduced also in the test setup.

Experimental characterization tests were carried out to assess the behavior of the bars against pull-off and the shear resistance of the GFRP connectors embedded in the mortar layer. The dimensioning of the steel bars was carried out so as to attain to a tensile resistance per unit of width slightly lower to that of the GFRP mesh. In such a way, the bars can yield and thus allow some rotation at the vault skewbacks, but providing some resistance against uplift and sliding. To let the bar yielding, an adequate bond length has to be guaranteed.

Both vaults with bricks arranged in a running bond pattern (two samples) and thin brick *in folio* vaults (four samples) were tested. The experimental behavior of the vaults samples was described and analyzed in the paper, also with indications of the crack patterns and of the capacity curves. The effectiveness of the reinforcement in improving both resistance and ductility was proved by comparison with the performances of unreinforced vaults: it clearly emerged that the presence of the GFRP mesh, which limits the opening of localized structural hinges at the haunches, delays the activation of the collapse mechanism, allowing instead a widespread mortar cracking. The steel bars, acting in tension, provided a partial rotational constraint at the skewback; some progressive extraction of the bars from the mortar was noted, however no visible uplift at the vault skewbacks occurred. Moreover, in the configuration reinforced at the extrados, the shear resistance against inward sliding resulted adequate. Differently, in the configuration reinforced at the intrados, the GFRP connectors allowed some inward sliding at skewbacks. Thus, the CRM reinforcement applied at the extrados resulted more effective than that applied at the intrados. However, in both cases, the performance improvement resulted significant (resistances about 4 to 6 times in running bond pattern vaults, 15 to 22 times in brick in folio vaults). Moreover, it was found that, in the "cracked" phase of the strengthened samples, the percentage of energy dissipation, in respect to the input energy, is significant (about 60%) and remained almost constant in all cycles up to collapse. It led to estimate a mean coefficient of equivalent viscous damping of 13%.

Further investigations can be performed for a higher optimization of the vault-to-wall connection strategy. For instance, the use of ribbed bars, instead of threaded ones, may improve the adhesion with the mortar layer and provide higher rotational constrain. Similarly, the use of more resistant or more numerous GFRP connectors may improve the shear resistance.

The obtained results will be useful for the refinement of the analytical method developed by the authors for the estimation of the reinforced vaults resistance [15], so to account for the influence of the connection with the abutments, useful for design purposes and reliable predictions. Moreover, the obtained results will allow to calibrate the hysteretic behavior of the strengthened sections and then to perform numerical analyses on different configurations, to deepen the study on the dissipative capacities of these elements.

## Acknowledgments

This paper is based on part of a research project financed by the composite engineering factory Fibre Net SpA, Pavia di Udine (UD), Italy. The authors wish also to thank Allen Dudine, Carlos Passerino and Andrea Miniussi for the useful help provided during the tests.